AAS 23-255

# DESIGN OF INITIAL GUESS LOW THRUST TRAJECTORIES USING CLOHESSY-WILTSHIRE EQUATIONS

Madhusudan Vijayakumar[*], William Skamser[†], Ossama Abdelkhalik[‡]

The commercial interest in producing low-cost space missions by exploiting the superior propellant management of low-thrust propulsion technology has become increasingly popular. Typical to such missions is the design of transfer trajectories between desired targets. This is a complex and computationally expensive process. Additionally, the optimal solvers used to generate these trajectories are extremely sensitive to initial guesses. One way to overcome this challenge is to use a reasonably approximate trajectory as an initial guess on optimal solvers. This paper presents a flexible approach to generating very low thrust trajectories. The initial guess is obtained from a flexible semi-analytic approach that can provide both planar and three-dimensional initial guess trajectories for various design scenarios like orbit raising, orbit insertion, phasing, and rendezvous. NASA's Evolutionary Mission Trajectory Generator (EMTG) and General Mission Analysis Tool (GMAT) are used as optimal solvers in this analysis. Numerical case studies are presented in this paper.

## INTRODUCTION

It is becoming increasingly popular to produce low-cost space missions by exploiting superior propellant management of low-thrust propulsion technology to perform rendezvous and proximity operations. Many agencies are considering the utilization of solar electric propulsion for satellite and debris removal missions in the near future. Additionally, this decade has seen a rise in the number of startups that plan to employ low-thrust propulsion systems to provide in-space servicing and refueling spacecraft. Likely, the space industry will eventually incorporate all-electric propulsion mission architecture for geocentric and cislunar missions. It is typical of such missions to design orbit maneuvers to transfer the spacecraft between desired targets. This design problem is usually solved using direct or indirect optimal control solvers. Direct solvers are preferred over indirect solvers due to their relatively larger radius of convergence. However, the initial guess required for direct solvers can be challenging to obtain due to the large number of design parameters.

Thus, the maneuver design of low-thrust missions is carried out in two stages. In the first stage, approximate initial guess trajectories are generated. This stage is often referred

[*]Ph.D. Candidate, Department of Aerospace Engineering, Iowa State University, USA; msudan@iastate.edu
[†]Undergraduate Student, Department of Aerospace Engineering, Iowa State University, USA; skamser@iastate.edu
[‡]Professor, Department of Aerospace Engineering, Iowa State University, USA; ossama@iastate.edu, AIAA Senior Member

to as the preliminary design. In the second stage, the initial guess trajectories are used to generate the high-fidelity optimal solution. This phase is also called the precise phase. The design space is rapidly explored during the preliminary stage to perform trade studies between feasibility and various mission parameters. This process is carried out for several weeks or even months before identifying trajectories that meet the mission requirements.

For a general low thrust trajectory optimization case, no closed-form analytic solution exists.[1] However, assuming some special condition, analytic approximation of orbit motion can be obtained,[2],[3].[4] On the numerical side, low-thrust trajectory design can be solved by treating them as a boundary value problem. A branch of trajectory design called the shape-based method was first introduced by Petropolous et al.[5] The shape-based method is a preliminary design where the trajectory shape is assumed to take up the form of a mathematical function. Following this work, a variety of shape-based methods have been developed. The exponential sinusoid[5] and the inverse polynomial[6] methods where the radial vector is expressed as a function of the transfer angle to generate initial guess orbit raising trajectories. The Bezier method[78] where the Bezier function describes the shape of the transfer trajectory. The shaping pseudoequinoctial,[9] and Finite Fourier Series[10] method extended the shape-based methods to generate trajectory solutions to a variety of problems like rendezvous, orbit insertion, etc., capable of handling thrust constraints. The finite Fourier series was later extensively extended to design planar,[10] three-dimensional,[11] suboptimal[1213] and trajectories in a multi-body environment[14].[15]

Several partial analytic solutions for special cases of orbit motion have been studied. One of the most important contributions to trajectory design was the derivation of the Clohessy-Wiltshire (CW) equations[16] in the 1960s that changed the game for trajectory design. This work provides the basis for analyzing relative spacecraft motion useful for proximity and rendezvous missions. For decades, these equations have been used as a reference model for designing guidance and navigation systems.[17] The proposed paper uses the CW equations to derive an analytic two-body low-thrust solution in the vicinity of a reference circular orbit. This analytic solution combines an iterative numerical solver to develop a fully automated shape-based initial trajectory generator. This method assumes that a constant magnitude thrust acceleration operates the spacecraft. The advantage of this method lies in its flexibility to handle various design scenarios like orbit raising, orbit insertion, orbit phasing, and rendezvous. Further, the reliability of the proposed method in generating quality initial guess trajectories is verified by interfacing it with NASA's Evolutionary Mission Trajectory Generator (EMTG) and General Mission Analysis Tool (GMAT) software tools.

**DYNAMIC MODEL**

The Clohessy-Wiltshire equations that describe the motion of a spacecraft relative to a reference orbit provide the dynamic framework for this work. The reference orbit is usually circular, and the corresponding coordinate frame of the relative motion is shown in Fig 1. The $X-Y-Z$ frame represents the inertial coordinate system, and $x-y-z$ represents the spacecraft's radial, along-track, and cross-track displacements relative to the reference

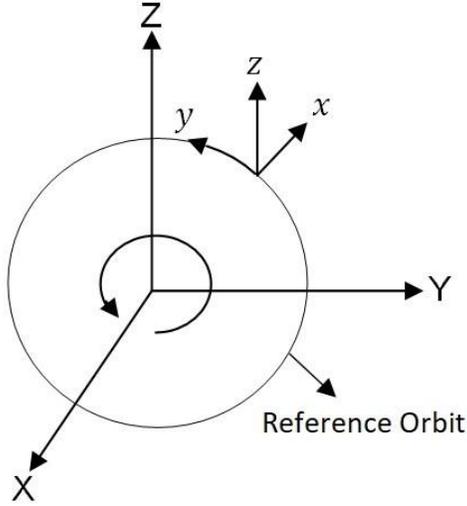 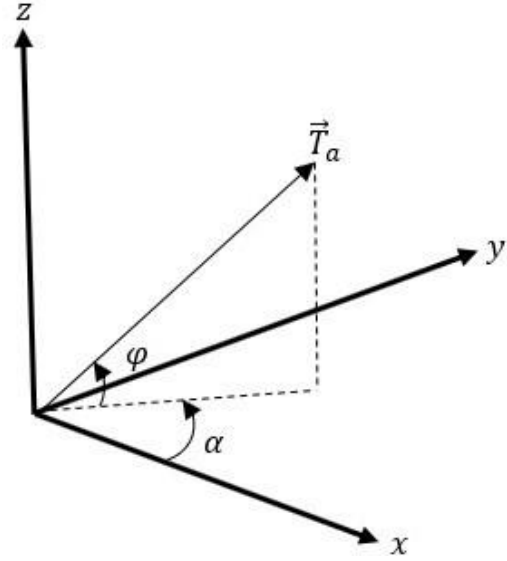

Figure 1. Coordinate Frames

Figure 2. Hills Frame

orbit. Consider the two-body Hill's equations as shown in Eq. (1).

$$\ddot{x} - 3xn^2 - 2\dot{y}n = a_x \quad (1a)$$
$$\ddot{y} + 2n\dot{x} = a_y \quad (1b)$$
$$\ddot{z} + zn^2 = a_z \quad (1c)$$

where $a_x$, $a_y$, and $a_z$ are the components of the trust acceleration represented in Hill's frame, and $n$ is the mean motion of the reference orbit. It is assumed that a constant magnitude thrust acceleration acts upon the spacecraft. Then the thrust acceleration components are given as:

$$a_x = a\cos(\phi)\cos(\alpha)$$
$$a_y = a\cos(\phi)\sin(\alpha) \quad (2)$$
$$a_z = a\sin(\alpha)$$

An alternate first-order matrix formulation of Hill's equations can be represented as follows:

$$\dot{\mathbf{x}} = \mathbf{f}(\mathbf{x}, \mathbf{u}) \quad (3)$$

where,

$$\mathbf{x} = \begin{bmatrix} x \\ y \\ z \\ \dot{x} \\ \dot{y} \\ \dot{z} \end{bmatrix}, \quad \mathbf{u} = \begin{bmatrix} \alpha \\ \phi \end{bmatrix}, \quad \mathbf{f}(\mathbf{x}, \mathbf{u}) = \begin{bmatrix} \dot{x} \\ \dot{y} \\ \dot{z} \\ 3xn^2 + 2\dot{y}n + a\cos(\phi)\cos(\alpha) \\ -2n\dot{x} + a\cos(\phi)\sin(\alpha) \\ -zn^2 + a\sin(\alpha) \end{bmatrix} \quad (4)$$

**SEMI-ANALYTIC SOLUTION**

This section describes a semi-analytic approach useful in generating initial guess low-thrust trajectories. An analytic solution describing spacecraft motion relative to a reference circular orbit is presented. This solution is combined with an iterative numerical search to achieve trajectory design.

**Analytic Solution**

An approximate analytic solution of the Clohessy-Wiltshire equations describing relative spacecraft motion is presented. The spacecraft is assumed to be acted upon by a low-thrust propulsion system providing constant magnitude thrust acceleration. Starting from the equations of motion described in Eq. (1) and using Laplace transforms, one can solve the set of linear differential equations and obtain time-domain solutions if the initial conditions are known. Noticing that the out-of-plane motion in Eq. (1c) is uncoupled from the in-plane motion, one can solve this separately. Converting Eq. (1c) to the Laplace domain, simplifying and using the initial conditions ($z_0$ and $\dot{z}_0$) one can obtain the time-domain solution of the out-of-plane motion. This can be expressed is as shown below:

$$z(t) = \frac{n^2 z_0 \cos(n\,t) - a \cos(n\,t) \sin(\beta) + n \dot{z}_0 \sin(n\,t)}{n^2} + \frac{a \sin(\beta)}{n^2} \tag{5}$$

$$\dot{z}(t) = \dot{z}_0 \cos(n\,t) - n z_0 \sin(n\,t) + \frac{a \sin(n\,t) \sin(\beta)}{n} \tag{6}$$

Given the initial conditions of in-plane states, ($x(0)$, $\dot{x}(0)$, $y(0)$, and $\dot{y}(0)$) one can use Eq. (1a) and (1b) and follow a similar approach to obtain an analytic solution describing spacecraft motion. This is shown below:

$$x(t) = f_x(x_0, \dot{x}_0, \alpha_0, \beta, n, k, t) \tag{7}$$

$$\dot{x}(t) = f_{\dot{x}}(x_0, \dot{x}_0, \alpha_0, \beta, n, k, t) \tag{8}$$

$$y(t) = f_y(x_0, \dot{x}_0, y_0, \alpha_0, \beta, k, t) \tag{9}$$

$$\dot{y}(t) = f_{\dot{y}}(x_0, \dot{x}_0, y_0, \alpha_0, \beta, k, t) \tag{10}$$

Detailed expressions for the on-track and along-track states ($x(t)$, $\dot{x}(t)$, $y(t)$, and $\dot{y}(t)$) represented in the Hill's frame is presented in the appendix.

**Trajectory Construction**

This section presents an iterative numerical process that utilizes the analytic solution for trajectory design. One should notice that the analytic approximation derived in the previous section is valid only for a small deviation from the reference orbit. This is because the fundamental model derived from these equations, namely the Clohessy-Wiltshire equations, is valid only for small deviations from the reference orbit. To maintain the predictability of the analytic solution, the entire trajectory is divided into 'm' segments and patched at points called the 'connection points.' Each segment is represented using an appropriate reference orbit, and the spacecraft states along each segment are approximated using the analytic solution. The segments are chosen so that the relative change in the orbital elements of adjacent segments is small. Figure 3 illustrates the segmentation process where the relative change in the orbital elements of adjacent segments is small. Dashed curves represent the reference orbits, blue circles represent the connection points, and solid curves represent the trajectory.

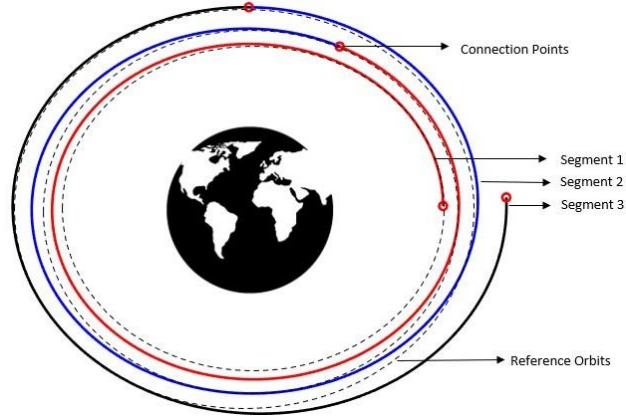

**Figure 3. Illustration for the Segmentation**

From Eq. (5)-Eq. (6) and Eq. (7)-(10) notice that the time of flight $t$, the thrust steering angles ($\alpha_0, \beta$) and the rate change of thrust acceleration direction $k$ are responsible for determining spacecraft motion. The desired transfer trajectories can be achieved by using a numerical solver to iterate over these parameters for the series of segments. The need for continuity conditions to maintain a smooth continuous trajectory is eliminated by enforcing the final states of the earlier segment to be the initial condition of the later segment.

## NUMERICAL EXAMPLES

This section the performance of the semi-analytic method as an initial guess generator for various trajectory design scenarios. NASA's Evolutionary Mission Trajectory Generator(EMTG) is chosen to be the optimal control tool used for testing. The low-thrust trajectory design problem on EMTG is set up as a finite burn low thrust (FBLT) mission. EMTG uses uniformly distributed discrete piece-wise control steps to solve this problem. To ensure constant magnitude thrust acceleration, control input on EMTG is set to the unit magnitude at all times. The process described in[18] is employed to input the initial guess generated by the semi-analytic approach on EMTG.

### Case 1: Earth to Mars Rendezvous Mission

A low-thrust Earth-Mars transfer problem is studied in this section. The spacecraft's departure date from Earth is arbitrarily chosen to be 20th July 2023. The orbital elements of Earth and Mars on this date are shown in Table 1. The DE430 Ephemeris model found on the Jet Propulsion Laboratory website[19] was used to obtain the planetary states. The spacecraft parameters are detailed in Table 2.

**Table 1. Case 1: Design Parameters**

| Parameter | Earth Orbit | Mars Orbit |
|---|---|---|
| Semi-major axis (AU) | 1 | 1.5236 |
| Eccentricity | 0.0167 | 0.0934 |
| Inclination (deg) | 0.00005 | 1.8506 |
| Longitude of Perihelion (deg) | 102.9471 | 336.0408 |
| Longitude of ascending node (deg) | -11.2606 | 49.5785 |
| True Anomaly (deg) | 194.72° | 201.99° |

**Table 2. Case 1: Spacecraft Parameters**

| Parameters | Value |
|---|---|
| Initial Mass (kg) | 1000 |
| Max Thrust (N) | 0.098 |
| Thrust-to-Weight ratio | 10^-6 |
| Isp | 2800 |

An Earth-Mars transfer solution was obtained by diving the trajectory design into 50 segments using the semi-analytic method. The trajectory solution and the history of the thrust steering angles are shown in Fig. 4 and Fig. 5 correspondingly. Using this solution

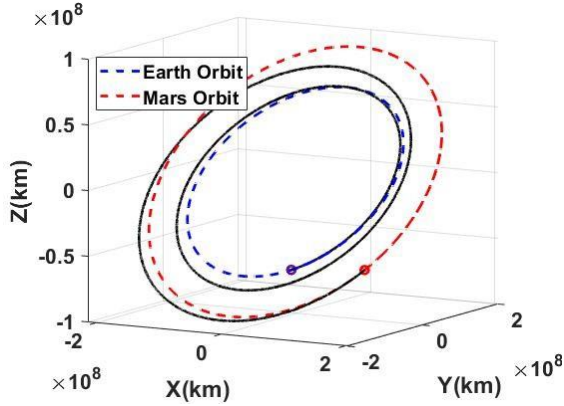
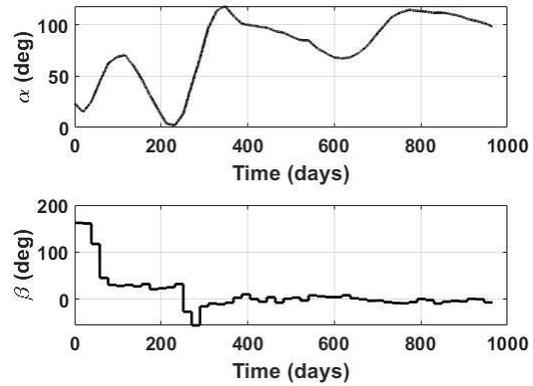

**Figure 4. Semi-Analytic: Earth-Mars Trajectory**    **Figure 5. Semi-Analytic: History of $\alpha$ and $\beta$**

as an initial guess and employing the Sparse Non-linear Optimizer (SNOPT) on EMTG, the converged solution is shown in Fig. 6. EMTG converged to the first feasible solution in 1.86 seconds after performing 3 major iterations. The corresponding thrust steering angles'

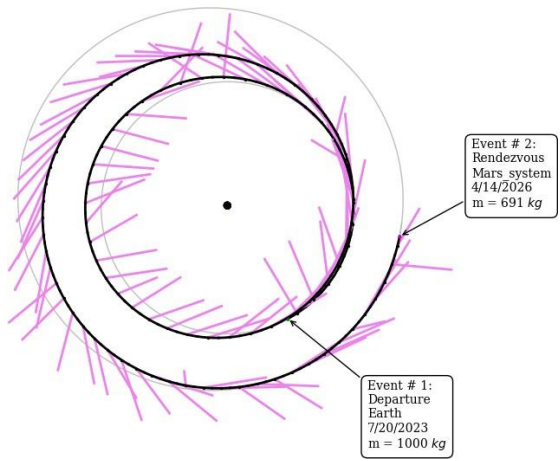

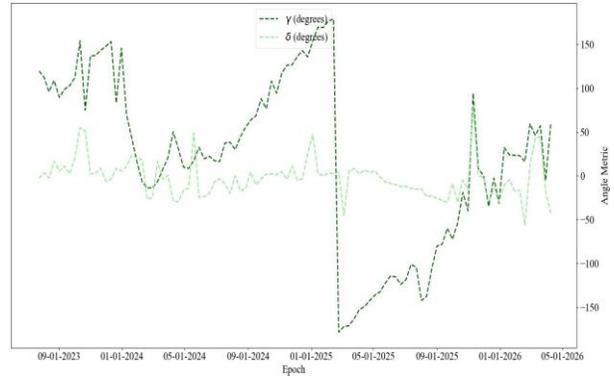

**Figure 6. EMTG: Earth-Mars Trajectory**

**Figure 7. EMTG: History of γ and δ**

history is shown in Fig. 7. Note that $\gamma$ and $\delta$ represent the in-plane and out-of-plane thrust steering angles in EMTG. The EMTG solution suggests that the spacecraft takes 999.97 days to transfer between Earth and Mars while consuming 308.57 kg of propellant. On the other hand, the time of flight and propellant mass were predicted to be 965.33 days and 261.87 days, corresponding to the semi-analytic solution. Similarly, the semi-analytic solution was used as an initial guess to solve an Earth-to-Mars rendezvous problem on GMAT. The only difference compared to the EMTG setup is the choice of design space. The time of flight window for this experiment was between 1 year to 10 years. The resulting trajectory consumed 287.8kg of propellant with a flight time of 933 days for this transfer. The trajectory solution and the history of the thrust steering angles are shown in Fig. 8 and Fig. 9 correspondingly. Table 3 draws a comparison between the initial guess trajectory

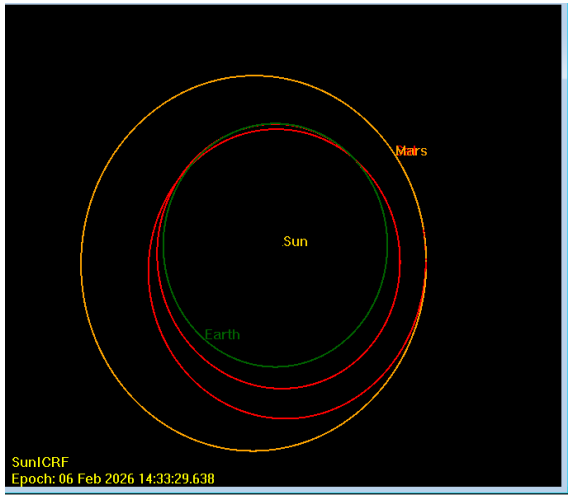

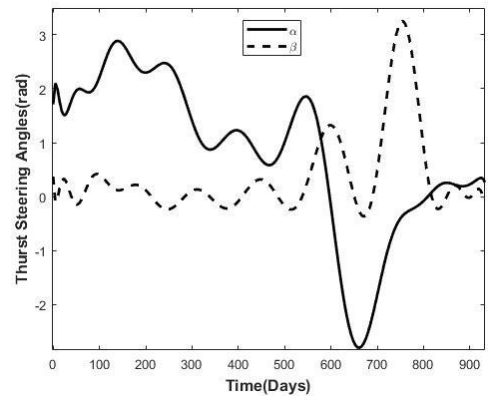

**Figure 8. GMAT: Earth-Mars Trajectory**

**Figure 9. GMAT: History of α and β**

solution obtained by the semi-analytic method and the converged trajectories on EMTG

and GMAT.

Table 3. Case 1: Comparison between the Semi-Analytic and EMTG Solutions

| Parameters | Semi-Analytic | EMTG | GMAT |
|---|---|---|---|
| No. of revolutions | 2 | 2 | 2 |
| Time of Flight (days) | 965.33 | 999.97 | 933 |
| Propellant Consumption (kg) | 261.87 | 308.57 | 287.8 |
| Computational Time (s) | 24.92 | 1.86 | - |

**Case 2: LEO to GEO Orbit Insertion Mission**

This section will study a low-thrust transfer between a circular low-earth orbit at an altitude of 300 km and the Geostationary Orbit. The parameters of the spacecraft are detailed in Table 4. The LEO-to-GEO transfer solution was obtained by dividing the trajectory design into 1440 segments. Since the number of segments is large, the design space of the optimization problem exploded. The trajectory was divided into multiple sections to overcome this challenge, and each section was solved sequentially. Sections consist of multiple segments and represent a portion of the trajectory. Sections were patched together to maintain continuity. The trajectory solution and the history of the thrust steering angles are shown in Figure 10 and 11. The spacecraft performs 369 revolutions around the Earth, taking approximately 62 days to reach GEO. The spacecraft consumes 14.01 *kg* of propellant, achieving a Δ*V* of 5.23 *km/s* for the transfer. The performance characteristics of the obtained solution are tabulated in Table 5. This LEO-to-GEO transfer solution was used as an initial guess on EMTG, with around 3690 control steps, i.e., 10 control steps per revolution. However, the experiment resulted in infeasible solutions. The experiment's outcome was not surprising, considering the massive dimensions of the design variables. The FBLT transcription on EMTG uses the control inputs specified as a three dimension vector (along the *x*,*y*, and *z* directions) at each control step. This results in 11,070 design variables to achieve the desired orbit transfer. For such a large design space, optimization becomes futile. On the other hand, decreasing the design variables would mean fewer control steps defining each revolution. What this means is reducing the flexibility of the optimizer in trying to achieve the desired target conditions. However, simply taking the semi-analytic solution and propagating it on NASA's General Mission Analysis Tool resulted in the trajectory shown in Figure 12. The final orbit achieved by the propagated trajectory reached is in the vicinity of the GEO with a semi-major axis of 42,047 km and an eccentricity of 0.001. This indicates the superior quality of the initial guess generated by the semi-analytic solution. Table 5 emphasizes the similarity between the two semi-analytic and their corresponding propagated trajectories on GMAT.

Table 4. Case 2: Spacecraft Parameters

| Parameters | Value |
|---|---|
| Initial Mass (kg) | 95 |
| Max Thrust Acceleration (m/s$^2$) | 9.81e-4 |
| Thrust-to-Weight ratio | 10^-4 |
| Isp | 3300 |

Table 5. Case 3: Performance Characteristics of the Semi-Analytic Solution

| Parameters | Value |
|---|---|
| No. of revolutions | 369 |
| Time of Flight (days) | 61.80 |
| $\Delta V$ (km/s) | 5.23 |
| Propellant Consumption (kg) | 14.01 |
| Computational Time (min) | 12 |

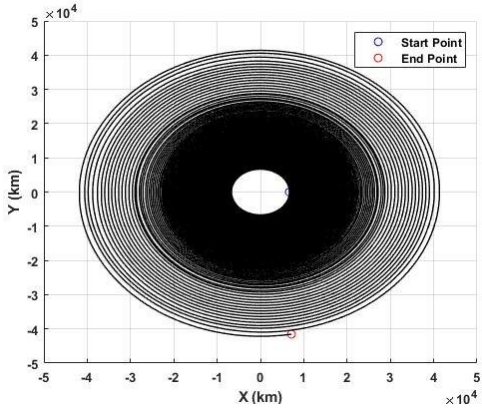 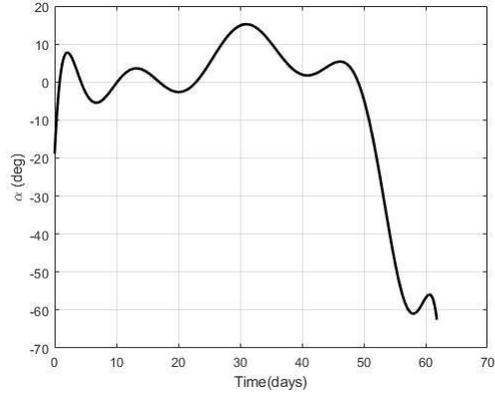

**Figure 10. Semi-Analytic:LEO-GEO Trajectory**  **Figure 11. Semi-Analytic: History of $\alpha$ and $\beta$**

**Table 6. Case 2: Comparison between the Semi-Analytic and GMAT Trajectories**

| Parameters | Semi-Analytic | GMAT |
|---|---|---|
| Final Semi-Major Axis (km) | 42165 | 42047 |
| Final Eccentricity | 0 | 0.001 |

**CONCLUSION**

The Semi-Analytic method successfully demonstrated its ability in generating reliable initial guess trajectories. Numerical results demonstrates the flexibility of the semi-analytic approach in generating three-dimensional rendezvous, orbit insertion, orbit phasing, and orbit raising trajectory solutions. Another contribution of this work is the ability in handling large number of orbital revolutions. This initial guess needs to be generated only once, and through successive seeding in an optimal solver higher-fidelity solutions can be obtained. The robustness of the semi-analytic solution combined with the streamlined numerical ap-proach makes it useful for automating trajectory design. The authors wish to extend this work to support cislunar missions in the future.

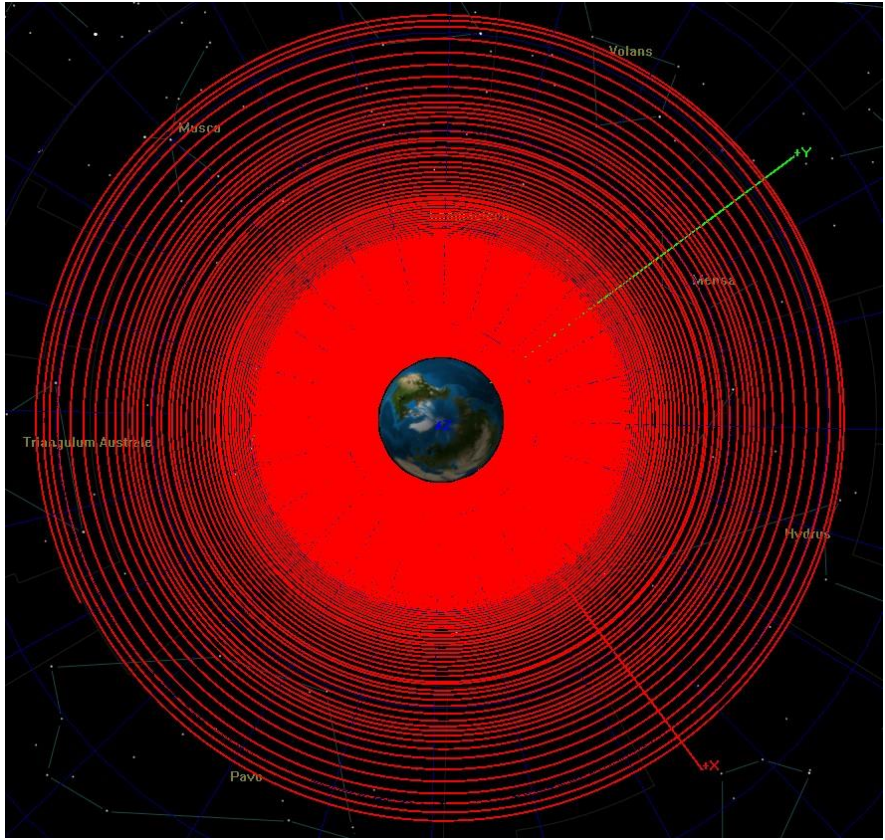

**Figure 12. GMAT: Trajectory Propagation LEO-GEO**

**ACKNOWLEDGMENT**

This paper is based upon work supported by NASA, Grant Number 80NSSC19K1642.

## Appendix - A

Analytic expression for spacecraft motion along the radial direction represented in Hill's frame.

$$x(t) = \frac{2 k^3 \dot{y}_0 - 4 k n^3 x_0 + 4 k^3 n x_0 - 2 k n^2 \dot{y}_0 - 2 k^3 \dot{y}_0 \cos(n t)}{k n (k^2 - n^2)}$$
$$+ \frac{k^3 \dot{x}_0 \sin(n t) - 2 a k^2 \cos(\alpha_0) \cos(\beta) + 2 a n^2 \cos(\alpha_0) \cos(\beta)}{k n (k^2 - n^2)}$$
$$+ \frac{3 k n^3 x_0 \cos(n t) - 3 k^3 n x_0 \cos(n t) + 2 k n^2 \dot{y}_0 \cos(n t)}{k n (k^2 - n^2)} \quad \text{(A.11)}$$
$$+ \frac{-k n^2 \dot{x}_0 \sin(n t) + 2 a k^2 \cos(n t) \cos(\alpha_0) \cos(\beta) - 2 a n^2 \cos(k t) \cos(\alpha_0) \cos(\beta)}{k n (k^2 - n^2)}$$
$$+ \frac{a k^2 \sin(n t) \cos(\beta) \sin(\alpha_0) - 2 a n^2 \sin(k t) \cos(\beta) \sin(\alpha_0) - a k n \cos(k t) \cos(\alpha_0) \cos(\beta)}{k n (k^2 - n^2)}$$
$$+ \frac{+a k n \cos(n t) \cos(\alpha_0) \cos(\beta) - a k n \sin(k t) \cos(\beta) \sin(\alpha_0) + 2 a k n \sin(n t) \cos(\beta) \sin(\alpha_0)}{k n (k^2 - n^2)}$$

$$\dot{x}(t) = \frac{k^2 \dot{x}_0 \cos(n t) - n^2 \dot{x}_0 \cos(n t) + 2 k^2 \dot{y}_0 \sin(n t)}{k^2 - n^2} + \frac{-3 n^3 x_0 \sin(n t) - 2 n^2 \dot{y}_0 \sin(n t) + 3 k^2 n x_0 \sin(n t)}{k^2 - n^2}$$
$$+ \frac{-a k \cos(k t) \cos(\beta) \sin(\alpha_0) + a k \sin(k t) \cos(\alpha_0) \cos(\beta) + a k \cos(n t) \cos(\beta) \sin(\alpha_0)}{k^2 - n^2} \quad \text{(A.12)}$$
$$+ \frac{-2 a k \sin(n t) \cos(\alpha_0) \cos(\beta) - 2 a n \cos(k t) \cos(\beta) \sin(\alpha_0) + 2 a n \sin(k t) \cos(\alpha_0) \cos(\beta)}{k^2 - n^2}$$
$$+ \frac{2 a n \cos(n t) \cos(\beta) \sin(\alpha_0) - a n \sin(n t) \cos(\alpha_0) \cos(\beta)}{k^2 - n^2}$$

## Appendix - B

Analytic expression for along-track spacecraft motion represented in Hill's frame.

$$y(t) = y_0 + \frac{t(3k^2 - 3n^2)(k\dot{y}_0 - a\cos(\alpha_0)\cos(\beta) + 2knx_0)}{k(n^2 - k^2)}$$

$$+ \frac{2k\sin(nt)(-3x_0 k^2 n - 2\dot{y}_0 k^2 + 2a\cos(\alpha_0)\cos(\beta)k + 3x_0 n^3 + 2\dot{y}_0 n^2 + a\cos(\alpha_0)\cos(\beta)n)}{k(n^3 - k^2)}$$

$$- \frac{2k\cos(nt)(\dot{x}_0 k^2 + a\cos(\beta)\sin(\alpha_0)k - \dot{x}_0 n^2 + 2a\cos(\beta)\sin(\alpha_0)n)}{k(n^3 - k^2)}$$

$$+ \frac{a\sin(\alpha_0 - kt)\cos(\beta)(k_0^2 - n^2)}{k^2(n^2 - k)} + \frac{2an\cos(kt)\cos(\beta)\sin(\alpha_0)(k + 2n)}{k^2(n^2 - k)}$$

$$- \frac{2an\sin(kt)\cos(\alpha_0)\cos(\beta)(k + 2n)}{k^2(n^2 - k)} - \frac{2\dot{x}_0 k^2 + 2a\cos(\beta)\sin(\alpha_0)k + 3an\cos(\beta)\sin(\alpha_0)}{k^2 n} \quad (B.13)$$

$$\dot{y}(t) = \dot{y}_0 + \frac{a\cos(\alpha_0 - kt)\cos(\beta)\sigma_1}{k(k^2 - n^2)}$$

$$- \frac{\cos(nt)(-6x_0 k^2 n - 4\dot{y}_0 k^2 + 4a\cos(\alpha_0)\cos(\beta)k + 6x_0 n^3 + 4\dot{y}_0 n^2 + 2a\cos(\alpha_0)\cos(\beta)n)}{k^2 - n^2}$$

$$- \frac{\sin(nt)(2\dot{x}_0 k^2 + 2a\cos(\beta)\sin(\alpha_0)k - 2\dot{x}_0 n^2 + 4a\cos(\beta)\sin(\alpha_0)n)}{k^2 - n^2} + \frac{2an\cos(kt)\cos(\alpha_0)\cos(\beta)(k + 2n)}{k(k^2 - n^2)}$$

$$+ \frac{2an\sin(kt)\cos(\beta)\sin(\alpha_0)(k + 2n)}{k(k^2 - n^2)} - \frac{4k\dot{y}_0 - 3a\cos(\alpha_0)\cos(\beta) + 6knx_0}{k} \quad (B.14)$$